\documentclass[a4paper,11pt,reqno]{amsart}
\usepackage{a4wide}

\usepackage[latin1]{inputenc}
\usepackage[T1]{fontenc}
\usepackage{amsmath}
\usepackage{amssymb}
\usepackage[all]{xy}

\newtheorem{theorem}{Theorem}

\begin{document}

\title{Fractional embeddings and stochastic time}

\author[J. Cresson]{J. Cresson}
\address[J. Cresson]{Laboratoire de Math\'ematiques appliqu\'ees de Pau, B\^atiment I.P.R.A., Universit\'e de Pau et des Pays de l'Adour, avenue de l'Universit\'e, BP 1155, 64013 Pau Cedex, France - Institut de M\'ecanique C\'eleste et de Calcul des \'Eph\'em\'erides, Observatoire de Paris, 77 avenue Denfert-Rochereau, 75014 Paris}
\email{jacky.cresson@univ-pau.fr}

\author[P. Inizan]{P. Inizan}
\address[P. Inizan]{Institut de M\'ecanique C\'eleste et de Calcul des \'Eph\'em\'erides, Observatoire de Paris, 77 avenue Denfert-Rochereau, 75014 Paris}
\email{pierre.inizan@imcce.fr}

\begin{abstract}                

As a model problem for the study of chaotic Hamiltonian systems, we look for the effects of a long-tail distribution of recurrence times on a fixed Hamiltonian dynamics.
We follow Stanislavsky's approach of Hamiltonian formalism for fractional systems.
We prove that his formalism can be retrieved from the fractional embedding theory.
We deduce that the fractional Hamiltonian systems of Stanislavsky stem from a particular least action principle, said causal. In this case, the fractional embedding becomes coherent.

\end{abstract}



\maketitle

\newcommand{\mc}[1]{\mathcal{#1}}
\newcommand{\Dc}[3]{{}_{#1} \mathcal{D}_{#2}^{#3} \, }
\newcommand{\Dcl}[0]{{}_{a} \mathcal{D}_{t}^{\alpha} \, }
\newcommand{\Dcr}[0]{{}_{t} \mathcal{D}_{b}^{\alpha} \, }
\newcommand{\Dpar}[2]{\dfrac{\partial #1}{\partial #2}}
\newcommand{\R}[0]{\mathbb{R}}
\newcommand{\psc}[2]{\langle #1 , #2 \rangle}
\newcommand{\Ker}[0]{\mathcal{K}er \,}
\newcommand{\Ima}[0]{\mathcal{I}m \,}


\section{Introduction}

Fractional calculus has been widely developped for a decade and its efficiency has already been proved in various topics such as continuum mechanics (see \cite{Carp_Main}), chemisty (see \cite{Bag_Tor_83}), transport theory (see \cite{Metzler_Klafter}), fractional diffusion (see \cite{Compte}), etc.

In \cite{Zasl_HCFD}, a link is drawn with chaotic Hamiltonian systems. Because of the appearance of fractal structures in phase spaces of nonhyperbolic Hamiltonian systems, fractional dynamics may arise in such systems. 
Zaslavsky then explained \cite[chap. 12-13]{Zasl_HCFD} that time takes on a fractal structure, meaning that it can be considered as a succession of specific temporal intervals. 
However, further investigations have to be carried out to understand and clarify the link between this peculiar temporal comportment and the fractional dynamics.

A contribution is done in \cite{Stan}. 
As a model problem for the effects of a given distribution of recurrence times on the underlying Hamiltonian dynamics, we use Stanislavsky's approach for his definition of an Hamiltonian formalism for fractional systems.
Indeed, this author looks for the effects induced by the assumption that the time variable is governed by a particular stochastic process on a given Hamiltonian dynamics. This kind of process contains notably the case of the algebraic decay of recurrence times that occurs in the study of chaotic Hamiltonian systems (see \cite{Zasl_HCFD}).
He proves, under strong assumptions, that the induced dynamics is fractional and that the structure of the new system looks like the classical Hamiltonian one. This allows him to give a definition of an Hamiltonian formalism for fractional systems.

However, an important property of Hamiltonian systems is that they can be obtained by a variational principle, called the \emph{Hamilton least action principle} (see \cite{Arnold}).
A natural question with respect to Stanislavsky's construction is to know if his definition of fractional Hamiltonian system can be derived from a variational principle.

In this paper, by using the \emph{fractional embedding} theory developped in \cite{Cresson}, we prove that Stanislavky's Hamiltonian formalism for fractional systems coincides with the fractional Hamiltonian formalism induced by the fractional embedding. In particular, this means that Stanis\-lavsky's fractional Hamiltonian systems can be obtained by a variational principle. Moreover, this fractional formalism is \emph{coherent}, meaning that there exists a commutative diagram for the obtention of the fractional equations.

In section \ref{Stan_form} we discuss Stanislavsky's formalism. 
Section \ref{Embedding} is devoted to the development of the fraction embedding theory using the Caputo derivatives. We obtain a causal and coherent embedding by restricting the set of variations underlying the fractional calculus of variations. We also prove that the fractional embedding of the usual Hamiltonian formalism resulting from the Lagrangian one is coherent. 
In section \ref{Compatibility}, we prove that the fractional Hamiltonian formalism stemming from the causal fractional embedding coincides with Stanislavsky's formalism.
We finally discuss open problems in section \ref{Conclusion}.


\section{Stanislavky's Hamiltonian formalism for fractional systems} \label{Stan_form}

\subsection{Definition of the internal time} \label{int_time}

Let $T_1, \, T_2, \, \ldots$ be nonnegative independant and identically distributed variables, with distribution $\rho$. We set $T(0) = 0$ and for $n \geq 1$, $T(n) = \sum_{i=1}^n T_i$. The $T_i$ represent random temporal intervals.
Let $\{N_t \}_{t \geq 0} = \max \{ n\geq 0 \, | \, T(n) \leq t \}$ be the associated counting process.
We suppose that there exists $ 0 < \alpha < 1$ such that
\begin{equation} \label{distrib}
\rho(t) \sim \frac{a}{t^{1+\alpha}}, \quad t \rightarrow \infty, \quad a > 0, \quad 0 < \alpha < 1. 
\end{equation} 
Therefore the variables $T_i$ belong to the strict domain of attraction of an $\alpha$-stable distribution.
Theorem 3.2 of \cite{Meer_Schef} implies:
\begin{theorem}
There exists a process $\{ S(t) \}_{t \geq 0}$ and a regularly varying function $b$ with index $\alpha$ such that
\begin{equation}
\{ b(c)^{-1} N_{ct} \}_{t \geq 0} \stackrel{FD}{\Longrightarrow} \{ S(t) \}_{t \geq 0}, \quad \text{as } c \rightarrow \infty \nonumber
\end{equation} 
where $\stackrel{FD}{\Longrightarrow}$ denotes convergence in distribution of all finite-dimensional marginal distributions.
\end{theorem}
The process $\{ S(t) \}_{t \geq 0}$ is a hitting-time process (see \cite{Meer_Schef}) and is also called a first-passage time.
From \cite{Bingham}, the distribution of $\{ S(t) \}_{t \geq 0}$, denoted $p_t$, verifies
 \begin{equation}
\mc{L}[p_t](v) = \mathbb{E}[e^{-v S(t)}] = E_\alpha(-v t^\alpha), \nonumber
\end{equation} 
where $\mc{L}$ is the Laplace transform and $E_\alpha$ is the one-parameter Mittag-Leffler function.
It follows that
\begin{equation} \label{Lapl}
\int_0^\infty e^{-w t} p_t(x) \,dt = w^{\alpha - 1} e^{-x w^\alpha} 
\end{equation}
The process $\{ S(t) \}_{t \geq 0}$ is increasing and may play the role of a stochastic time, which is called \emph{internal time} in \cite{Stan}. The distribution $p_t(\tau)$ represents the probability to be at the internal time $\tau$ on the real time $t$.
Using this new time, Stanislavsky studies Hamiltonian systems which evolve according to $S(t)$.

\subsection{Fractional Hamiltonian equations}

We consider an Hamiltonian system, with Hamiltonian $H(x,p)$, and associated canonical equations
\begin{equation} \label{eq_can}
\begin{array}{rcl} 
 \dfrac{d}{dt} x(t) & = & \partial_2 H (x(t),p(t)), \\
 \dfrac{d}{dt} p(t) & = & - \partial_1 H (x(t),p(t)). 
\end{array}
\end{equation} 

If $t$ is replaced by $S(t)$, how is the dynamics modified? 
To answer this question, Stanislavsky introduces new variables $x_\alpha$ and $p_\alpha$ defined by
\begin{equation} \label{xpa}
\begin{array}{rcl} 
x_\alpha (t) = \mathbb{E}[x(S(t))] & = & \int_0^\infty p_t(\tau) x(\tau) d\tau, \\
p_\alpha (t) = \mathbb{E}[p(S(t))] & = & \int_0^\infty p_t(\tau) p(\tau) d\tau. 
\end{array}
\end{equation} 

Furthermore, he assumes that 
\begin{equation} \label{hypothese_H}
\begin{array}{rcl} 
\partial_1 H (x_\alpha(t), p_\alpha(t)) & = & \int_0^\infty p_t(\tau) \partial_1 H (x(\tau),p(\tau)) d\tau, \\
\partial_2 H (x_\alpha(t), p_\alpha(t)) & = & \int_0^\infty p_t(\tau) \partial_2 H (x(\tau),p(\tau)) d\tau, 
\end{array}
\end{equation} 

which leads to
\begin{theorem} \label{Dyn_Stan}
Let $(x,p)$ be a solution of \eqref{eq_can}. Then condition \eqref{hypothese_H} is verified if and only if $(x_\alpha,p_\alpha)$ defined by \eqref{xpa} verifies
\begin{equation} \label{Stan_eq}
\begin{array}{rcl} 
\Dc{0}{t}{\alpha} x_\alpha (t) & = & \partial_2 H (x_\alpha(t), p_\alpha(t)), \\
\Dc{0}{t}{\alpha} p_\alpha (t) & = & - \partial_1 H (x_\alpha(t), p_\alpha(t)), \\
\end{array}
\end{equation} 

where $\Dcl$ is the left Caputo derivative defined by
\begin{equation}
\Dcl f(t) = \frac{1}{\Gamma(1-\alpha)} \int_a^t (t-\tau)^{-\alpha} f'(\tau) \, d\tau. \nonumber
\end{equation} 

\end{theorem}

\begin{proof}
As $x$ verifies \eqref{eq_can}, we have 
\begin{equation}
\int_0^\infty p_t(\tau) \partial_2 H(x(\tau),p(\tau)) \, d\tau = \int_0^\infty p_t(\tau) \frac{d}{d\tau} x(\tau) d\tau \nonumber
\end{equation} 
The Laplace transform of this expression gives
\begin{eqnarray}
\mc{L} \left[ \int_0^\infty p_t(\tau) \partial_2 H(x,p) \, d\tau \right](s)  = \int_0^\infty \mc{L}[p_t](s) \frac{d}{d\tau} x(\tau) d\tau \nonumber \\
\qquad  = s^{\alpha-1} \int_0^\infty e^{-\tau s^\alpha} \frac{d}{d\tau} x(\tau) d\tau \; \text{ from }  \eqref{Lapl} \nonumber \\
  = s^{2\alpha-1} \mc{L}[x](s^\alpha) - s^{\alpha-1} x(0). \qquad \qquad \!  \nonumber
\end{eqnarray} 
Given that $\mc{L}[x_\alpha](s) = s^{\alpha-1} \mc{L}[x](s^\alpha)$, we have
\begin{equation}
\mc{L} \left[ \int_0^\infty p_t(\tau) \partial_2 H(x(\tau),p(\tau)) \, d\tau \right](s) = \mc{L} \left[ \Dc{0}{t}{\alpha} x_\alpha \right](s) \nonumber.
\end{equation} 
By taking the Laplace image of this relation, we obtain 
\begin{equation}
	\Dc{0}{t}{\alpha} x_\alpha(t) = \int_0^\infty p_t(\tau) \partial_2 H(x(\tau),p(\tau)) \, d\tau. \nonumber
\end{equation}
In a similar way, we also have
\begin{equation}
	\Dc{0}{t}{\alpha} p_\alpha(t) = - \int_0^\infty p_t(\tau) \partial_1 H(x(\tau),p(\tau)) \, d\tau, \nonumber
\end{equation}
and the equivalence follows.
\end{proof}

For a presentation of the fractional calculus and its applications, see \cite{Samko} and \cite{Oldham}.
Hence, we will say that a fractional system of the form
\begin{equation} 
\begin{array}{rcl} 
\Dc{0}{t}{\alpha} x (t) & = & f_1 (x(t), p(t)), \nonumber \\
\Dc{0}{t}{\alpha} p (t) & = & f_2 (x(t), p(t)), \nonumber \\
\end{array}
\end{equation} 

is Hamiltonian in the sense of Stanislavsky if there exists a function $H(x,p)$ such that 
\begin{equation} 
\begin{array}{rcl} 
f_1(x,p)  & = & \partial_2 H(x,p), \nonumber \\
f_2(x,p) & = & - \partial_1 H(x,p). \nonumber \\
\end{array}
\end{equation}

We show that the fractional derivative $\Dc{0}{t}{\alpha}$ appears as a natural consequence of the structure of the internal time $S(t)$. The fractional exponent $\alpha$ is exactly determined by the behaviour \eqref{distrib} of long time intervals. We note that if we had $\alpha \geq 1$ in \eqref{distrib}, the $\alpha$-stable distribution would be the Gaussian one, we would have $p_t(\tau) = \delta_\tau(t)$ and then $S(t) \equiv t$. In this case, internal time and real time would be the same. Consequently, for $\alpha \geq 1$, the associated derivative is the classical one.


\section{Fractional embedding of Lagrangian and Hamiltonian systems} \label{Embedding}

One important property of classical Hamiltonian systems is that they are solutions of a variational principle, called the \emph{Hamilton least action principle} (see \cite{Arnold}). A natural question is to know if the fractional Hamiltonian systems defined by Stanislavky can be derived from a variational principle.

Fractional Euler-Lagrange and Hamilton equations has been first derived in \cite{Riewe_96}, in order to include frictional forces into a variational principle. In \cite{Agrawal}, a fractional Euler-Lagrange equation is obtained using a fractional least action principle. This formalism includes left and right fractional derivatives. The related Hamilton equations are derived in \cite{Baleanu_Agrawal}. However, their equations are different from those obtained by Stanislavsky.

Using the fractional embedding theory developped in \cite{Cresson}, we prove that Stanislavsky Hamiltonian formalism stems from a fractional variational principle, called \emph{causal}, and moreover that this construction is coherent.

We sum up here the general ideas of the fractional embedding theory for the Caputo derivative. 
Similarly to the left one, the right Caputo derivative is defined by
\begin{equation}
\Dc{t}{b}{\alpha} \, f(t) = \frac{-1}{\Gamma(1-\alpha)}  \int_t^{b} (\tau- t)^{-\alpha} f'(\tau) \, d\tau. \nonumber
\end{equation}

The left fractional integral is defined by
\begin{equation}
\Dc{a}{t}{- \alpha} f(t) = \frac{1}{\Gamma(\alpha)} \int_a^t (t-\tau)^{\alpha-1} f(\tau) \, d\tau, \nonumber
\end{equation}

and the right one by
\begin{equation}
\Dc{t}{b}{- \alpha} f(t) = \frac{1}{\Gamma(\alpha)} \int_t^b (\tau-t)^{\alpha-1} f(\tau) \, d\tau. \nonumber
\end{equation}

\subsection{Fractional embedding of differential operators}

Let $\textbf{f} = (f_1, \ldots, f_p)$ and $\textbf{g} = (g_1, \ldots, g_p)$ be two $p$-uplets of smooth functions $\R^{k+2} \longrightarrow \R^l$. Let $a,b \in \R$ with $a < b$.
We denote $\mc{O} (\textbf{f},\textbf{g})$ the differential operator defined by
\begin{equation} \label{Ofg}
\mc{O} (\textbf{f},\textbf{g}) (x)(t) = \sum_{i=0}^p ( f_i \cdot \frac{d^i}{dt^i} g_i ) ( x(t), \ldots, \frac{d^k}{dt^k} x(t), t ), 
\end{equation} 
where, for any functions $f$ and $g$, $(f \cdot g) (t) = f(t) \cdot g(t)$, where $\cdot$ means a product component by component.

The fractional embedding of $\mc{O} (\textbf{f},\textbf{g})$, denoted $\mc{E}_\alpha(\mc{O} (\textbf{f},\textbf{g}))$, is defined by
\begin{equation} \label{emb}
\mc{E}_\alpha(\mc{O} (\textbf{f},\textbf{g}))(x)(t) = \! \!\sum_{i=0}^p ( f_i \cdot ( \Dcl \!)^i g_i ) ( x(t), .., ( \Dcl \!)^k x(t), t ) \nonumber
\end{equation}   

We define the ordinary differential equation associated to $\mc{O} (\textbf{f},\textbf{g})$ by
\begin{equation} \label{ODE}
\mc{O} (\textbf{f},\textbf{g})(x) = 0. 
\end{equation} 

The fractional embedding $\mc{E}_\alpha(\mc{O} (\textbf{f},\textbf{g}))$ of \eqref{ODE} is defined by
\begin{equation}
\mc{E}_\alpha(\mc{O} (\textbf{f},\textbf{g})) (x) = 0. \nonumber
\end{equation}

\subsection{Lagrangian systems}

Now we consider a Lagrangian system, with smooth Lagrangian $L(x,v,u)$ and $u \in [a,b]$. The Lagrangian $L$ can naturally lead to a differential operator of the form \eqref{Ofg}:
\begin{equation}
\mc{O} (1,L)(x)(t) =  L ( x(t), \frac{d}{dt} x(t), t ). \nonumber
\end{equation} 
Now we identify $L$ and $\mc{O}(1,L)$.
The fractional embedding of $L$, $\mc{E}_\alpha(L)$, is hence given by
\begin{equation}
\mc{E}_\alpha(L)(x)(t) = L( x(t), \Dcl x(t), t). \nonumber
\end{equation}   

In Lagrangian mechanics, the action and its minima play a central role. For any mapping $g$, the action of $g$, denoted $\mc{A}(g)$ is defined by
\begin{equation}
\mc{A}(g)(x) = \int_a^b g(x)(t) \, dt. \nonumber
\end{equation} 

For example, with the identification $L \equiv \mc{O}(1,L)$, the action of $L$ is given by
\begin{equation}
\mc{A}(L)(x) = \int_a^b L( x(t), \dfrac{d}{dt} x(t), t ) \, dt, \nonumber
\end{equation} 

and concerning the fractional embedding of $L$, the associated action is
\begin{equation}
\mc{A}(\mc{E}_\alpha(L))(x) = \int_a^b L \left( x(t), \Dcl x(t), t \right) \, dt. \nonumber
\end{equation} 

The extremum of the action of a Lagrangian $L$ provides the equation of motion associated:
\begin{theorem}
The action $\mc{A}(L)$ is extremal in $x$ if and only if $x$ satisfies the Euler-Lagrange equation, given by
\begin{equation} \label{EL_clas}
\partial_1 L ( x(t), \dfrac{d}{dt} x(t), t ) - \dfrac{d}{dt}  \partial_2 L ( x(t), \dfrac{d}{dt} x(t), t ) = 0.
\end{equation}
This equation is denoted $EL(L)$.
\end{theorem}

This procedure should not be modified with fractional derivatives. Indeed, the strict definition of the Lagrangian $L$ does not involve any temporal derivative. The dynamics is afterwards fixed with the choice of the derivative $\mc{D}$ and the relation $v(t) = \mc{D} x(t)$.
The variational principle providing the Euler-Lagrange equation uses a integration by parts, which remains in the fractional case:
\begin{multline} \label{IPP}
\int_a^b \left[ \Dc{a}{t}{\alpha} f(t) \right] g(t) dt = \int_a^b f(t) \left[ \Dc{b}{t}{\alpha} g(t) \right] dt \\
+ g(b) \Dc{a}{b}{-(1-\alpha)} f(b) - f(a) \Dc{a}{b}{-(1-\alpha)} g(a). \nonumber
\end{multline}

We introduce the space of variations 
\begin{equation}
V_\alpha = \{ h \in C^1([a,b]) \: | \: \Dc{a}{b}{-(1-\alpha)} h(a) = h(b) = 0 \}. \nonumber
\end{equation} 
For $h \in V_\alpha$, we have
\begin{multline}
\mc{A}(\mc{E}_\alpha(L))(x+h) = \mc{A}(\mc{E}_\alpha(L))(x) + \\
\int_a^b [\partial_1 L + \Dcr \partial_2 L](x(t), \Dcl x(t), t) \, h(t) \, dt + o(h), \nonumber
\end{multline} 

which implies that the differential of $\mc{A}(\mc{E}_\alpha(L))$ in $x$ is given, for any $h \in V_\alpha$, by 
\begin{align}
d \mc{A}(\mc{E}_\alpha(L))(x,h) \! & = \! \! \int_a^b [\partial_1 L \! + \! \Dcr \partial_2 L](x, \Dcl x, t) \, h(t) \, dt, \nonumber \\
				   & = \psc{[\partial_1 L \! + \! \Dcr \partial_2 L](x(\cdot), \Dcl x(\cdot), \cdot)}{h}, \nonumber
\end{align} 

where $\psc{f}{g} = \int_a^b f(t) g(t) \, dt$ is a scalar product defined on $C^1([a,b])$.

If $E \subset V_\alpha$, we will say that $\mc{A}(\mc{E}_\alpha(L))$ is $E$-extremal in $x$ if for all $h \in E$, $d \mc{A}(\mc{E}_\alpha(L))(x,h) = 0$.

So we obtain a first Euler-Lagrange equation:
\begin{theorem}
$\mc{A}(\mc{E}_\alpha(L))$ is $V_\alpha$-extremal in $x$ if and only if $x$ verifies 
\begin{equation} \label{EL1}
\partial_1 L(x(t), \Dcl x(t), t)\! + \! \Dcr \partial_2 L(x(t), \Dcl x(t), t) = 0.
\end{equation} 
\end{theorem}

\begin{proof}
$\mc{A}(\mc{E}_\alpha(L))$ is $V_\alpha$-extremal in $x$ if and only if for all $h \in V_\alpha$, we have
\begin{equation*}
\psc{[\partial_1 L \! + \! \Dcr \partial_2 L](x(\cdot), \Dcl x(\cdot), \cdot)}{h} = 0.
\end{equation*} 

This is equivalent to $[\partial_1 L \! + \! \Dcr \partial_2 L](x(\cdot), \Dcl x(\cdot)) \in V_\alpha^\perp$. We conclude by noticing that $V_\alpha^\perp = \overline{V_\alpha}^\perp = \{ 0 \}$, where $\overline{V_\alpha}$ is the adherence of $V_\alpha$ in $C^1([a,b])$, equal to $C^1([a,b])$ entirely.
\end{proof}

Equation \eqref{EL1} will be called \emph{general fractional Euler-Lagrange equation} and will be denoted $EL_g(\mc{E}_\alpha(L))$.
Contrary to \eqref{EL_clas}, two operators are involved here. We will now discuss the problematic presence of $\Dcr$.

\subsection{Coherence and causality}

Because of the simultaneous presence of the two derivatives, the position of $x$ at time $t$ depends on its past positions, through $\Dcl$, but also on its future ones, through $\Dcr$.
The principle of causality is here violated, which seems crippling from a physical point of view. 
Moreover, we note that \eqref{EL_clas} can be written in the form \eqref{ODE}, with $\textbf{f}=(1,1)$ and $\textbf{g}=(\partial_1 L, - \partial_2 L)$. The fractional embedding $\mc{E}_\alpha(EL(L))$ of \eqref{EL_clas} is therefore
\begin{equation} 
\partial_1 L(x(t), \Dcl x(t), t)\! - \! \Dcl \partial_2 L(x(t), \Dcl x(t), t) = 0, \nonumber
\end{equation} 
which shows that $EL_g(\mc{E}_\alpha(L)) \not\equiv \mc{E}_\alpha(EL(L))$: fractional embedding and least action principle are not commutative.
So we obtain two procedures providing different fractional equations, which seems also unsatisfactory.
 We are facing a Cornelian choice: shall we preserve causality or the least action principle?  
A possible way to solve this problem is to restrict the space of variations. We note $\tilde{V}_\alpha = \{ h \in V_\alpha \, | \, \Dcl h = - \Dcr h \}$ and $K_\alpha = \Dcl + \Dcr$, defined on $C^1([a,b])$. For any $f, \, g \in V_\alpha$, $\psc{K_\alpha \, f}{g} = \psc{f}{K_\alpha \, g}$. We show that $K_\alpha$ is essentially self-adjoint and we obtain a new Euler-Lagrange equation:
\begin{theorem}
$\mc{A}(\mc{E}_\alpha(L))$ is $\tilde{V}_\alpha$-extremal in $x$ if and only if there exists a function $g$ such that $x$ verifies 
\begin{equation} \label{EL2}
\partial_1 L(x(t), \Dcl x(t), t)\! - \! \Dcl \partial_2 L(x(t), \Dcl x(t), t) = K_\alpha \, g. \nonumber
\end{equation} 
\end{theorem}

\begin{proof}
$\mc{A}(\mc{E}_\alpha(L))$ is $\tilde{V}_\alpha$-extremal in $x$ if and only if $[\partial_1 L \! + \! \Dcr \partial_2 L](x(\cdot), \Dcl x(\cdot)) \in \tilde{V}_\alpha^\perp$. Given that $\tilde{V}_\alpha^\perp = (\Ker K_\alpha)^\perp = \Ima K_\alpha$, $\mc{A}(\mc{E}_\alpha(L))$ is extremal if and only if there exists $\tilde{g}$ such that 
\begin{equation} 
\partial_1 L(x(t), \Dcl x(t), t)\! + \! \Dcr \partial_2 L(x(t), \Dcl x(t), t) = K_\alpha \, \tilde{g}. \nonumber
\end{equation} 
We conclude by setting $g(t) = \tilde{g}(t) + \partial_2 L(x(t), \Dcl x(t), t)$. 
\end{proof}

Restricting of the space of variations breaks the unicity of the solution. However, among those solutions, there is a single one which remains causal (without the operator $\Dcr$), for $g=0$:
\begin{equation} \label{EL3}
\partial_1 L(x(t), \Dcl x(t), t)\! - \! \Dcl \partial_2 L(x(t), \Dcl x(t), t) = 0. 
\end{equation} 
Equation \eqref{EL3} will be called \emph{causal fractional Euler-Lagrange equation}, and will be denoted $EL_c(\mc{E}_\alpha(L))$.

Now causality is respected and we have $EL_c(\mc{E}_\alpha(L)) \equiv \mc{E}_\alpha(EL(L))$. In this case, the fractional embedding is told \emph{coherent}, in the sense that the following diagram commutes:

\begin{eqnarray}
\xymatrix{
 L \ar[d]_{\mbox{(causal) LAP}} \ar[r]^{\mc{E}_\alpha} & \mc{E}_\alpha(L)
 \ar[d]^{\mbox{causal FLAP}}       \\
 \left( \partial_1 L - \frac{d}{dt} \partial_2 L \right) = 0 \ar[r]_{\mc{E}_\alpha}  & \left( \partial_1 L - \Dcl \partial_2 L \right) = 0  } \nonumber
\end{eqnarray}

where (F)LAP states for "(fractional) least action principle". As $\Dc{a}{t}{1} = -\Dc{t}{b}{1} = \dfrac{d}{dt}$, we can say that the least action is also causal in the classical case.

However, in the fractional case, the physical meaning of $\tilde{V}_\alpha$ is not clear, but it might be related to a reversible dynamics of the variations.  
Furthermore, this underlines the significant role of variations in the global dynamics.

\subsection{Fractional Hamiltonian systems based on fractional Lagrangian ones} \label{Subsec_Leg}

There exists a natural derivation of an Hamiltonian system from a Lagrangian system based on the Legendre transformation. 
We consider an autonomous Lagrangian system, with Lagrangian $L(x,v)$, and we suppose that
\begin{equation} \label{Legendre}
\forall x, \; v \mapsto \partial_2 L(x,v) \text{ is bijective}.
\end{equation}
The (static) momentum associated to the (static) variable $x$ is $p=\partial_2 L(x,v)$. So there exists a mapping $f$ named Legendre transformation such that $v = f(x,p)$. The Hamiltonian $H$ associated to $L$ is defined by 
\begin{equation} \label{Ham}
H(x,p) = p f(x,p) - L(x, f(x,p)). 
\end{equation} 
It implies $\partial_1 H(x,p) = - \partial_1 L(x, f(x,p))$ and $\partial_2 H(x,p) = f(x,p)$. 

Let introduce the function
\begin{equation}
F_{LH}(x,p,v,w) = \begin{pmatrix}
p - \partial_2 L(x,v) \nonumber \\
\partial_1 H(x,p) + \partial_1 L(x, f(x,p)) \nonumber \\
\partial_2 H(x,p) - f(x,p) \nonumber
\end{pmatrix}
\end{equation}

The link between Lagrangian and Hamiltonian formalisms is done through the equation 
\begin{equation} \label{FLH}
F_{LH}(x,p,v,w) = 0.
\end{equation}

The momentum $p$ induces a function $p(t) \! = \! \partial_2 L(x(t),v(t))$, which can be considered as the dynamical momentum.

For the classical dynamics, \eqref{FLH} becomes
\begin{equation}
F_{LH}\left( x(t),p(t),\dfrac{d}{dt}x(t),\dfrac{d}{dt}p(t) \right) = 0, \nonumber
\end{equation}
i.e.
\begin{eqnarray}
p(t) & = & \partial_2 L \left( x(t),\dfrac{d}{dt}x(t) \right), \nonumber \\
\partial_1 H(x(t),p(t)) & = & - \partial_1 L \left( x(t), \dfrac{d}{dt}x(t) \right), \nonumber \\
\partial_2 H(x(t),p(t)) & = & \dfrac{d}{dt}x(t). \nonumber
\end{eqnarray}

Moreover, if $x(t)$ is solution of the Euler-Lagrange equation \eqref{EL_clas}, we obtain the canonical equations 
\begin{equation}
\begin{array}{rcl} 
 \dfrac{d}{dt} x(t) & = & \partial_2 H (x(t),p(t)), \nonumber \\
 \dfrac{d}{dt} p(t) & = & - \partial_1 H (x(t),p(t)). \nonumber
\end{array}
\end{equation} 

For the fractional case, the fractional embedding \eqref{emb} of \eqref{FLH} is given by
\begin{eqnarray}
p(t) & = & \partial_2 L \left( x(t), \Dcl x(t) \right), \nonumber \\
\partial_1 H(x(t),p(t)) & = & - \partial_1 L \left( x(t), \Dcl x(t) \right), \nonumber \\
\partial_2 H(x(t),p(t)) & = & \Dcl x(t). \nonumber
\end{eqnarray}

Then the following result states:
\begin{theorem} \label{Leg_H}
If $x(t)$ is solution of the causal fractional Euler-Lagrange equation \eqref{EL3}, we have
\begin{equation}
\begin{array}{rcl} 
 \Dcl x(t) & = & \partial_2 H (x(t),p(t)), \nonumber \\
 \Dcl p(t) & = & - \partial_1 H (x(t),p(t)). \nonumber
\end{array}
\end{equation} 
\end{theorem}

These are the equations describing the dynamics of a fractional Hamiltonian system derived from a Lagrangian formalism. But Hamiltonian systems can also be considered directly as it will be seen.  

\subsection{Embedded Hamiltonian systems}

Now we consider an Hamiltonian system as defined in section \ref{Stan_form}, with an Hamiltonian $H(x,p)$ and with equations \eqref{eq_can} associated. The fractional embedding \eqref{emb} of \eqref{eq_can} is
\begin{equation}
\begin{array}{rcl} \label{Ham_emb}
 \Dcl x(t) & = & \partial_2 H (x(t),p(t)), \\
 \Dcl p(t) & = & - \partial_1 H (x(t),p(t)). 
\end{array}
\end{equation}

Furthermore, by indroducing the function 
\begin{equation}
L_H(x,p,v,w) = pv - H(x,p), \nonumber
\end{equation}

we see that classical Hamiltonian systems are critical points of the action of $L_H$ defined by
\begin{equation}
\mc{A}(L_H)(x,p) = \int_a^b L_H \left( x(t),p(t), \dfrac{d}{dt} x(t), \dfrac{d}{dt} p(t) \right) \, dt. \nonumber
\end{equation}
 
In the fractional case, the action becomes
\begin{equation}
\mc{A}(\mc{E}_\alpha(L_H))(x,p) \! = \! \! \int_a^b \! L_H \left( x(t),p(t), \Dcl x(t), \Dcl p(t) \right) dt. \nonumber
\end{equation}

Using the causal fractional Euler-Lagrange equation for $L_H$, we obtain

\begin{theorem}[Hamiltonian coherence]
Let $H$ be an Hamiltonian function. The solutions $(x(t),p(t))$ of the fractional system \eqref{Ham_emb} coincide with causal critical points of the action $\mc{A}(\mc{E}_\alpha(L_H))$.
More precisely, the following diagram commutes:

\begin{eqnarray}
\xymatrix{
 & L_H \ar[d]_{\mbox{(causal) LAP}} \ar[r]^{\mc{E}_\alpha} & \mc{E}_\alpha(L_H)
 \ar[d]^{\mbox{causal FLAP}}       \\
 & {\left\{ \begin{array}{rcl} 
 \tfrac{d}{dt} x(t) \! & \! = \! & \! \partial_2 H \\
 \tfrac{d}{dt} p(t) \! & \! = \! & \! - \partial_1 H
\end{array}\right .} \ar[r]_{\mc{E}_\alpha}  & 
{\left \{ \begin{array}{rcl} 
 \Dcl x(t) \! & \! = \! & \! \partial_2 H \\
 \Dcl p(t) \! & \! = \! & \! - \partial_1 H 
 \end{array}\right .} } \nonumber
\end{eqnarray} 
\end{theorem}

\begin{proof}
The causal fractional Euler-Lagrange equation for $L_H$ is
\begin{eqnarray}
- \partial_1 H(x(t),p(t)) - \Dcl p(t) & = & 0, \nonumber \\
\Dcl x(t) - \partial_2 H(x(t),p(t)) & = & 0, \nonumber
\end{eqnarray}

which is exactly \eqref{Ham_emb}.

\end{proof}

So we have coherence between the directly embedded equations and the equations obtained by a variational principle. But we have also coherence between this section and the previous one, i.e. between the fractional Hamiltonian systems resulting from Lagrangian ones and the embedded Hamiltonian systems.

In other words, the equivalent approaches for Hamiltonian systems in the classical case remain equivalent in the fractional case if we use causal variational principles. 

Now we will discuss the link between this formalism and Stanislavsky's one.
 

\section{Compatibility between the two formalisms} \label{Compatibility}

Condition \eqref{hypothese_H} means that the partial derivatives of $H$ commutes with $\mathbb{E}[\cdot(S(t))]$. This condition could just seem of technical order and could appear as unrelated to the real dynamics. However, by using the fractional embedding, we can precise the underlying dynamical link, in the case of natural Lagrangian systems. 

We consider a natural Lagrangian system, i.e. with a Lagrangian $L$ of the form $L(x,v) = \dfrac{1}{2}m v^2 - U(x)$, and the Hamiltonian $H(x,p)$ derived as in section \ref{Subsec_Leg}. So we have $H(x,p) = \dfrac{1}{2m} p^2 + U(x)$, with $p = \partial_2 L(x,v) = mv$. We suppose that $(x,p)$ is solution of the classical Hamiltonian equations \eqref{eq_can}. We define the associated variables $x_\alpha$ and $p_\alpha$ by \eqref{xpa}. 
  
\begin{theorem}
 If $x_\alpha$ is solution of the causal fractional Euler-Lagrange equation \eqref{EL3} associated to $L$, then condition \eqref{hypothese_H} is verified.
\end{theorem}

\begin{proof}
We set $\tilde{p}_\alpha(t) = \partial_2 L(x_\alpha(t), \Dc{0}{t}{\alpha} x_\alpha(t))$, i.e. $ \tilde{p}_\alpha(t) = m \: \Dc{0}{t}{\alpha} x_\alpha(t)$. Then, from theorem \eqref{Leg_H}, $(x_\alpha,\tilde{p}_\alpha)$ is solution of 
\begin{equation} \label{Comp_pf}
\begin{array}{rcl} 
 \Dc{0}{t}{\alpha} x_\alpha(t) & = & \partial_2 H (x_\alpha(t),\tilde{p}_\alpha(t)), \nonumber \\
 \Dc{0}{t}{\alpha} \tilde{p}_\alpha(t) & = & - \partial_1 H (x_\alpha(t),\tilde{p}_\alpha(t)). \nonumber
\end{array}
\end{equation} 

Moreover, we have 
\begin{multline}
\tilde{p}_\alpha(t) = m \: \Dc{0}{t}{\alpha} \! \int_0^\infty \! \! p_t(\tau) x(\tau) \, d\tau = m \int_0^\infty \! \! p_t(\tau) \dfrac{d}{d\tau} x(\tau) \, d\tau  \\
= \int_0^\infty p_t(\tau) m v(\tau) \, d\tau = \int_0^\infty p_t(\tau) p(\tau) \, d\tau = p_\alpha(t).
\end{multline}

So we can replace $\tilde{p}_\alpha$ by $p_\alpha$ in \eqref{Comp_pf}, which concludes the proof.

\end{proof}


\section{Conclusion} \label{Conclusion}

If we consider the temporal evolution variable of a Lagrangian system as a succession of random intervals, and if their density has a power-law tail, then the dynamics of this system is fractional. The associated equations can be determined through a fractional embedding, based on a least action principle. In order to obtain causal and coherent equations, it is necessary to restrict the space of variations. This condition might be seen as a way to cancel the finalist aspect of the least action principle. 
Even if it is still unclear, this model of time could notably be appropriated for the description of some chaotic Hamiltonian dynamics. 
Some numerical experiments show that distributions of Poincaré recurrence times may possess a power-law tail (see \cite[chap. 11]{Zasl_HCFD}, \cite{Zasl_Fan}, \cite{Cristadoro}). Consequently, the time may be decomposed into a succession of recurrence times. For long time scale dynamics, the number of intervals is great and the new characteristic time clock may become $S(t)$. This new time takes into account the peculiar structure of the recurrence times: if the power-law exponent $\alpha$ verifies $0 < \alpha < 1$, the long time scale dynamics becomes fractional with the same exponent $\alpha$. This idea of stacked dynamics based on two time scales could be linked with \cite{Hilfer_FFD}, where close results are obtained. 
However, because of the Kac lemma (\cite{Kac}), which states that the mean recurrence time is finite, condition \eqref{distrib} may be valid only locally, near some island boundaries, called sticky zones. Further investigations have to be carried on to clarify this point.


\bibliographystyle{plain}
\bibliography{Biblio_Cresson_Inizan}             

\begin{thebibliography}{10}

\bibitem{Agrawal}
O.P. Agrawal.
\newblock Formulation of {E}uler-{L}agrange equations for fractional
  variational problems.
\newblock {\em J. Math. Anal. Appl.}, 272:368--379, 2002.

\bibitem{Arnold}
V.I. Arnold.
\newblock {\em Mathematical methods of classical mechanics}.
\newblock Springer, 1989.

\bibitem{Bag_Tor_83}
R.L. Bagley and P.J. Torvik.
\newblock A theoritical basis for the application of fractional calculus in
  viscoelasticity.
\newblock {\em Journal of Rheology}, 27:201--210, 1983.

\bibitem{Baleanu_Agrawal}
D.~Baleanu and O.P. Agrawal.
\newblock Fractional {H}amilton formalism within {C}aputo's derivative.
\newblock {\em Czechoslovak Journal of Physics}, 56(10-11), 2006.

\bibitem{Bingham}
N.H. Bingham.
\newblock Limit theorems for occupation times of {M}arkov processes.
\newblock {\em Z. Warscheinlichkeitsth.}, 17:1--22, 1971.

\bibitem{Carp_Main}
A.~Carpinteri and F.~Mainardi.
\newblock {\em Fractals and Fractional Calculus in Continuum Mechanics}.
\newblock Springer-Verlag, New York, 1997.

\bibitem{Compte}
A.~Compte.
\newblock Stochastic foundations of fractional dynamics.
\newblock {\em Physical Review E}, 53(4):4191--4193, 1996.

\bibitem{Cresson}
J.~Cresson.
\newblock Fractional embedding of differential operators and {L}agrangian
  systems.
\newblock {\em J. Math. Phys.}, 48(3):033504, 2007.

\bibitem{Cristadoro}
G.~Cristadoro and R.~Ketzmerick.
\newblock Universality of algebraic decays in {H}amiltonian systems.
\newblock {\em Physical Review Letters}, 100(18):184101(4), 2008.

\bibitem{Zasl_Fan}
R.~Fan and G.M. Zaslavsky.
\newblock Pseudochaotic dynamics near global periodicity.
\newblock {\em Comm. Non. Sci. Num. Sim.}, 12:1038--1052, 2007.

\bibitem{Hilfer_FFD}
R.~Hilfer.
\newblock Foundations of fractional dynamics.
\newblock {\em Fractals}, 3(3):549--556, 1995.

\bibitem{Kac}
M.~Kac.
\newblock {\em Probability and Related Topics in Physical Sciences}.
\newblock Interscience, New York, 1957.

\bibitem{Meer_Schef}
M.M. Meerschaert and H.-P. Scheffler.
\newblock Limit theorems for continuous-time random walks with infinite mean
  waiting times.
\newblock {\em J. Appl. Prob.}, 41:623--638, 2004.

\bibitem{Metzler_Klafter}
R.~Metzler and J.~Klafter.
\newblock The random walk's guide to anomalous diffusion: a fractional dynamics
  approach.
\newblock {\em Pysics Reports}, 339:1--77, 2000.

\bibitem{Oldham}
K.B. Oldham and J.~Spanier.
\newblock {\em The Fractional Calculus}.
\newblock Academic Press, New York and London, 1974.

\bibitem{Riewe_96}
F.~Riewe.
\newblock Nonconservative {L}agrangian and {H}amiltonian mechanics.
\newblock {\em Physical Review E}, 53(2):1890, 1996.

\bibitem{Samko}
S.G. Samko, Kilbas A.A., and O.I. Marichev.
\newblock {\em Fractional integrals ans derivatives: theory and applications}.
\newblock Gordon and Breach, New York, 1993.

\bibitem{Stan}
A.A. Stanislavsky.
\newblock Hamiltonian formalism of fractional systems.
\newblock {\em Eur. Phys. J. B}, 49:93--101, 2006.

\bibitem{Zasl_HCFD}
G.M. Zaslavsky.
\newblock {\em Hamiltonian Chaos \& Fractional Dynamics}.
\newblock Oxford University Press, Oxford, 2005.

\end{thebibliography}

\end{document}